\documentclass[11pt]{article}
\usepackage{amsfonts}

\usepackage{graphics}
\usepackage{indentfirst}
\usepackage{cite}
\usepackage{latexsym}
\usepackage{amsmath}
\usepackage{amssymb}
\usepackage[dvips]{epsfig}
\usepackage{amscd}

\newtheorem{theorem}{Theorem}[section]
\newtheorem{remark}{Remark}[section]

\newtheorem{lemma}[theorem]{Lemma}

\newtheorem{corollary}[theorem]{Corollary}

\newcommand{\n}{\rho}

\newcommand{\ti}{\tilde}

\def\pf{{\it Proof.}  }

\newcommand{\pa}{\partial}

\newcommand{\bi}{\bibitem}

\newcommand{\bt}{\begin{theorem}}
\newcommand{\bl}{\begin{lemma}}
\newcommand{\el}{\end{lemma}}
\newcommand{\et}{\end{theorem}}
\newcommand{\ga}{\gamma}

\newcommand{\de}{\delta}
\newcommand{\ve}{\varepsilon}
\newcommand{\la}{\label}
\newcommand{\si}{\sigma}

\newcommand{\bn}{\begin{eqnarray}}
\newcommand{\en}{\end{eqnarray}}
\newcommand{\bnn}{\begin{eqnarray*}}
\newcommand{\enn}{\end{eqnarray*}}

\newcommand{\bnnn}{\begin{eqnarray*}}
\newcommand{\ennn}{\end{eqnarray*}}
\newcommand{\ben}{\begin{enumerate}}
\newcommand{\een}{\end{enumerate}}

\newcommand{\ba}{\begin{aligned}}
\newcommand{\ea}{\end{aligned}}
\newcommand{\be}{\begin{equation}}
\newcommand{\ee}{\end{equation}}

\def\norm[#1]#2{\|#2\|_{#1}}

\def\xix{\int_{-\infty}^{+\infty}}

\def\xiT{\int_0^T}

\def\xR{\mathbb{R}}

\def\xl{\left}
\def\xr{\right}

\makeatletter      
\@addtoreset{equation}{section}
\makeatother       

\title{Global Well-Posedness and Large-Time  Behavior of   1D   Compressible Navier-Stokes System with Density-Depending Viscosity and Vacuum in Unbounded Domains\thanks{K.  Li is partially supported by Undergraduate Research Fund of BNU 2017-150; B. L\"u is supported by NNSFC (Nos. 11601218 \& 11771382), Science and Technology Project  of Jiangxi Provincial Education Department (No. GJJ160719).}}

\author{  Kexin Li\thanks{
 School of Mathematical Sciences,
 Beijing Normal University,
  Beijing 100875,  P. R. China ({\tt kexinli98@163.com}) }   \quad  Boqiang L\"u\thanks{College of Mathematics and Information Science, Nanchang Hangkong University, Nanchang 330063, P. R. China({\tt lvbq86@163.com}). }
  \quad Yixuan Wang\thanks{ Institute of Applied Mathematics, AMSS,
Chinese Academy of Sciences, Beijing 100190,  P. R. China({\tt wangyixuan\_14@163.com}).
 }
 }
 \date{}

\begin{document}
\maketitle

\begin{abstract} We consider the  Cauchy problem for one-dimensional (1D) barotropic compressible Navier-Stokes equations with density-dependent viscosity and large external force. Under a general  assumption on the density-dependent  viscosity,  we prove that the Cauchy problem admits a unique  global strong (classical) solution   for the  large initial data with vacuum. Moreover,    the density is proved to be bounded from  above time-independently. As a consequence,      we   obtain the large time behavior of the solution without external forces.  \end{abstract}



\section{Introduction and main results}

The motion of a one-dimensional viscous compressible barotropic fluid is governed by the following compressible Navier-Stokes equations:
\be\la{R1d}
\begin{cases} \rho_t + (\rho u)_x = 0,\\
 (\rho u)_t + (\rho u^2)_x + [P(\n)]_x =[\mu(\n)u_x]_x+\n f.
\end{cases}
\ee
 Here, $t\ge 0$ is time, $x\in \mathbb{R}=(-\infty,\infty)$ is the spatial coordinate,
 and $\n(x,t)\ge 0$,  $u(x,t)$  and
$P(\rho) = A\rho^{\gamma} ( A>0,\, \ga>1)$
are the fluid density, velocity and pressure, respectively. Without loss of generality, it is assumed that $A= 1$. The viscosity   $\mu(\cdot)$  is a function of the density $\n$.  The external force $f=f(x)$ is a  known function. 
We look for the solutions $(\n,u)$ to the Cauchy problem for \eqref{R1d} with the  initial condition
\be\la{R1d1}
(\n,u)|_{t=0}=(\n_0(x),u_0(x)),~~~\mbox{for}~x\in  \xR,
\ee
and   the following  far field behavior
\be\la{R1d2}
u(x,t)\rightarrow0,~~ \n(x,t)\rightarrow\tilde{\n}>0,~~~\mbox{as}~|x|\rightarrow \infty,
\ee
where $\tilde\n$ is a given positive constant.

There is huge literature on the studies of the global existence and large time behavior of solutions to the 1D compressible Navier-Stokes equations.
For constant viscosity $\mu$ and the initial data away from vacuum, the problems are addressed by Kanel \cite{Kaz} for sufficiently smooth data, and
by Serre \cite{serre1,serre2} and Hoff \cite{Hof} for discontinuous initial data.  On the other hand, if  $\mu$ depends on $\n$ and admits a positive constant lower bound, the  global well-posedness and large time behavior of solutions away from vacuum were discussed in \cite{zl89,v1989,is2002} and the references therein. 
However, when it comes to the case that vacuum is allowed initially,
as emphasized in many papers related
to compressible fluid dynamics \cite{L1,Fe,coi1,xin98,hl,hlx1}, the possible presence of vacuum is one of the major
difficulties in discussing the well-posedness of solutions  to the
compressible Navier-Stokes equations.
In the presence of vacuum, Ding-Wen-Zhu \cite{wen2011}  considered  the global existence of classical  solutions to  1D compressible Navier-Stokes equations in bounded domains,  provided that $\mu \in C^2[0,\infty)$ satisfies
\be\la{wen}
 0<\bar{\mu}\le \mu(\n)\le C(1+P(\n)).
\ee Recently, for general $\mu,$ we  \cite{lww} establish not only the global existence but also the large-time behavior for   classical   solutions containing  vacuum  to the initial boundary value problem for  1D compressible Navier-Stokes equations.
  For the  Cauchy problem \eqref{R1d}--\eqref{R1d2} without external force $(f=0),$ Ye \cite{ye15} studies the global classical large solutions  under the following restriction on $\mu(\n)$:
\be\la{ye}
\mu(\n)=1+\n^\beta,~~~0\le\beta<\gamma.
\ee
However,  both the uniform upper bound of density
and the large time behavior of  solutions are not obtained in \cite{ye15}.

In this paper, for more general density-dependent  viscosity (see \eqref{Rn3}) and large external force,  we will derive the uniform upper bound of density and thus prove  the global well-posedness of strong (classical) large solutions containing  vacuum   to  the  Cauchy problem \eqref{R1d}--\eqref{R1d2}.
Before stating the main results, we first explain the notations and
conventions used throughout this paper.
 We set$$D_t\triangleq \frac{\pa}{\pa t} +u  \frac{\pa }{\pa x},~~~\dot v\triangleq v_t+uv_x.$$
For   $1\le r\le \infty$ and $k\ge 1$, we adopt the following simplified notations:
$$L^r=L^r(\mathbb{R}),\quad W^{k,r}  = W^{k,r}(\mathbb{R}) , \quad H^k = W^{k,2}(\mathbb{R}).$$

The first result   concerns the global existence   of  strong solutions to the Cauchy problem \eqref{R1d}--\eqref{R1d2}.
\begin{theorem}\la{Rth1}   Suppose that $\int_{-\infty}^x f(y)dy\in H^2$ and  the viscosity   $\mu(\n)\in C^1[0,\infty)$   satisfies
\be\la{Rn3}
 0<\bar{\mu}\le \mu(\n)\le \lambda_0\int_1^\n \mu(s)ds+\lambda_1,
\ee for some  constants $\bar{\mu}>0,~\lambda_0\ge 0,$ and $\lambda_1>0$. Let  the initial data $(\n_0,u_0)$ satisfy
\be\la{R1d3}
\n_0\ge 0,~ \n_0- \tilde{\n} \in H^1, ~u_0\in    H^1. \ee
Then, there exists a unique global strong
solution $(\n,u)$ to the Cauchy problem \eqref{R1d}--\eqref{R1d2} satisfying  for any $0< T<\infty,$
   \be
   \la{R1d5}\begin{cases}
    \n-\tilde{\n} \in C([0,T]; H^1), \quad \n_t\in L^\infty(0,T;L^2),\\
 u\in L^\infty(0,T; H^1) \cap L^2(0,T;H^2),\\
 \sqrt{t}u \in L^\infty(0,T; H^2),~  \sqrt{t} u_t\in L^2(0,T;H^1).\end{cases} \ee  Moreover, the density remains uniformly bounded for all time, that is
 \be\la{R1d6} \sup_{0\le t<\infty} \|\n(\cdot,t)\|_{L^\infty}<\infty.
\ee
\end{theorem}

   Similar to    \cite[Theorem 1.2]{lww}, one can obtain directly the following result which shows that the strong solutions obtained by Theorem \ref{Rth1} become classical provided initial data $(\n_0,u_0)$ satisfy some additional  compatibility  conditions.

\begin{theorem}\la{Rth11}   In addition to the conditions of Theorem \ref{Rth1},  suppose that       $\mu(\n)\in C^2[0,\infty)$  and  that the initial data $(\n_0,u_0)$ satisfy
\be\la{xR1d3}
 (\n_0- \tilde{\n},~P(\n_0)-P(\tilde{\n}))\in H^2, ~u_0\in    H^2, \ee
and the   compatibility condition:
 \be\la{R1d4}
[\mu(\n_0)u_{0x}]_x-[P(\n_0)]_x=\sqrt{\n_0}g(x),~~~x\in \xR,\ee
for a given function $g\in L^2$. Then, the strong solution obtained in Theorem \ref{Rth1} becomes classical and satisfies  for any $0< \tau< T<\infty,$
   \be
   \la{xR1d5}\begin{cases}
   (\n-\tilde{\n},~P(\n)-P(\tilde{\n})) \in C([0,T]; H^2), \\
 u\in C([0,T]; H^2) \cap L^2(\tau,T;H^3),\\
   u_t\in L^2(0,T;H^1),~ \sqrt{t}u \in L^\infty(0,T; H^3),\\
    \sqrt{t}u_{t}\in L^\infty(0,T; H^1)\cap L^2(0,T; H^2),~ \sqrt{t}\sqrt{\n}u_{tt}\in L^2(0,T; L^2).\end{cases} \ee
\end{theorem}

When there is no external force, that is $f\equiv 0$ in \eqref{R1d}, we can obtain  the large time behavior of the strong solutions to  the Cauchy problem \eqref{R1d}--\eqref{R1d2}.

\begin{theorem}\la{Rth2} In addition to the conditions of   Theorem \ref{Rth1}, suppose that $f\equiv 0.$
Then the following large time behavior holds for strong solution $(\n, u)$ to
the Cauchy problem \eqref{R1d}--\eqref{R1d2} obtained by Theorem \ref{Rth1}:
\be\la{R1d7}
\lim_{t\rightarrow\infty}\left( \|\n -\tilde \n\|_{L^p}+ \|u_x\|_{L^2\cap L^p}\right)=0,~~~~ \forall~p>2.
\ee
 Moreover, if  there exists some point $x_0\in (-\infty,\infty)$ such that $\n_0(x_0)=0,$  the spatial gradient of the density the unique strong solution $(\n, u)$ to  the  problem  \eqref{R1d}--\eqref{R1d2}  has to blow up as $t\rightarrow \infty$ in the following sense,
\be\label{bp}
\lim_{t\rightarrow \infty}\|\n_x(\cdot,t)\|_{L^2}=\infty .
\ee
\end{theorem}

A few remarks are in order:

\begin{remark}\la{Rre1}
It should be noted here that we obtain a completely  new  uniform upper bound of density \eqref{R1d6}, which is in sharp contrast to \cite{wen2011,ye15} where the upper bound of the  density are time-dependent. This is crucial  for studying the large time behavior of the solutions. 
\end{remark}

\begin{remark} \la{Rre2}
We want to point out  that our restriction  on $\mu(\n)$ in \eqref{Rn3} is much more general in comparison  with those of \cite{ye15,wen2011} (see \eqref{ye} and \eqref{wen}). Here, we will list some special cases satisfying \eqref{Rn3} as follows:

$\bullet$  $\mu(\n)=\hat{\mu}$ is a positive constant. Let $\lambda_0=0$ and $\lambda_1=\hat{\mu}.$ Then
$$\lambda_0\int_1^\n \mu(s)ds+\lambda_1=\hat{\mu}=\mu(\n).$$

$\bullet$ $\mu(\n)=1+\n^a$ for any $a\ge 0$. Choosing $\lambda_0=1+a$ and $\lambda_1=4+a$ leads to
$$\lambda_0\int_1^\n \mu(s)ds+\lambda_1=(1+a)\n+\n^{1+a}+2\ge 1+\n^a=\mu(\n).$$

$\bullet$ $\mu(\n)=e^\n$. Letting $\lambda_0=1$ and $\lambda_1=e $ yields
$$\lambda_0\int_1^\n \mu(s)ds+\lambda_1=e^\n=\mu(\n).$$
Hence, it is clear that  the density-dependent viscosity $\mu(\n)$ in \cite{ye15,wen2011} are all included in our results.
\end{remark}

\begin{remark} \la{xre4} To obtain the global  existence  of strong solutions in Theorem \ref{Rth1}, we do not need the additional  compatibility condition \eqref{R1d4}, which is required for discussing the global classical solutions in Theorem \ref{Rth11}. Theorems \ref{Rth1}--\ref{Rth11} show that how the compatibility condition \eqref{R1d4} plays its role in studying the well-posedness of solutions with initial vacuum.
\end{remark}

\begin{remark}\la{re1}
It should be noted here that the solution  $(\n,u)$  obtained in Theorem \ref{Rth1} is actually a classical one.
Indeed,   the Sobolev embedding theorem  together with the regularities of  $(\n,u)$  in \eqref{xR1d5}  shows that
\be \label{ccc}(\n,\,P,\, u) \in C([0 ,T];C^{1+\frac{1}{2}}),~~ \n_t \in C([0 ,T];C^{\frac{1}{2}}).\ee
Furthermore, one can deduce from \eqref{xR1d5} that  for any $0<\tau<T$,
$$u \in L^\infty(\tau,T; H^3),~~~u_t \in L^\infty(\tau,T; H^1)\cap L^2(\tau,T; H^2),$$
which yields that for $1<r<2$,
\be \label{ccc2}u \in   C([\tau,T]; C^{2,r}[0,1]),~~~u_t \in   C([\tau,T]; C^r[0,1]).\ee
Hence, one can deduce from \eqref{ccc}--\eqref{ccc2}   that  $(\n,u)$   is  the classical solution to   \eqref{R1d}--\eqref{R1d2}.
\end{remark}

We now make some comments on the analysis of this paper.
Since  the local well-posedness result of classical solutions away from vacuum is well-known  (see  \cite{coi1} or Lemma  \ref{Rpro1} below),  to   extend   local solutions  to be a global one which will eventually  contain  vacuum after letting the  lower bound of initial density go to zero,  we need some
global a priori estimates which are independent of the lower bound of density. As mentioned in many papers (see \cite{lij06jpma,lzz}),  the main difficulty comes from     initial vacuum, density-dependent viscosity,  large external force, and the unboundedness of the domain.  Motivated by our previous result \cite{lww}, we find that
the key issue is to derive both the time-independent upper bound of the  density and
the time-depending derivative ones of the solutions. Motivated by \cite{ka1},
we   localize the problem on bounded domains and use the method develop in our previous work \cite{lww} to get the upper bound of the density independently of $x$ and thus the uniform bound of density on the whole $\mathbb{R}$ due to the  arbitrariness   of  $x$ (see Lemma \ref{Rla3.2}).
Furthermore, with the help of far field behavior of $\n$, we can  bound the $L^2$-norm of $u$ in terms of $\|\n^{1/2}u\|_{L^2}$ and $\|u_x\|_{L^2}$ (see \eqref{R3.5}).   Following the methods used in  \cite{lzz,lx1},  we use the material derivative  $\dot u$
instead of the usually $u_t$ and succeed in obtaining  the  derivative estimates on $(\n, u).$
  With the desired global a priori estimates independent of the lower bound of density at hand, we can thus prove that  Theorems \ref{Rth1} and \ref{Rth11} hold  for   any time $T>0$ which completes the proof of  Theorems \ref{Rth1} and \ref{Rth11} (see Section 3). Finally, for the case without  external force ($f\equiv 0$), we can establish the time-independent lower order  estimates (see \eqref{xRgj3.5} and \eqref{xRgj3.6}).  Using the methods due to \cite{hlx1,lzz} and the key  time-independent a priori estimates,  we prove  Theorem \ref{Rth2}  in Section 4.

The rest of the paper is organized as follows.   Section 2   is devoted to derive the necessary
a priori estimates on smooth solutions.
Theorems \ref{Rth1}  and \ref{Rth2} are  proved in Sections 3 and 4, respectively.

 \section{A priori estimates}

 First,  for the initial density    strictly away from vacuum, we state the following  local well-posedness theory of strong solutions,   whose proof can be obtained by similar arguments as in \cite{coi1,M1}.
 \begin{lemma}\la{Rpro1} Let $\mu(\n)\in C^2[0,\infty).$  Assume that    $\tilde f\in H^2$ and the initial data $(\n_0,u_0)$ satisfies
   \bnn\la{Rlocal1}
 (\n_0- \tilde{\n},~P(\n_0)-P(\tilde{\n}),u_0)\in    H^2,~\inf_{x\in\xR}\n_0(x)>0,
\enn
where  $\tilde\n>0$ is a given constant.
Then, there exists a small time $T_0>0$  such that the Cauchy problem \eqref{R1d}--\eqref{R1d2}   has a unique strong solution $(\n,u)$ on $\xR\times(0, T_0]$ satisfying
   \bnn
   \la{Rpre3}\begin{cases}
\n-\tilde{\n},~P(\n)-P(\tilde{\n})\in C([0,T_0]; H^2), \\
 u\in C([0,T_0]; H^2) \cap L^2( 0,T_0 ; H^3),\\
  \n_t\in C([0,T_0]; H^1),~~u_t\in C([0,T_0]; L^2)\cap L^2( 0,T_0 ; H^1).
 \end{cases} \enn
\end{lemma}

In this section, we will establish some  a priori bounds
for strong solutions $(\n,u)$  to the Cauchy problem \eqref{R1d}--\eqref{R1d2} on $\xR\times [0,T]$ whose existence is guaranteed   by Lemma \ref{Rpro1}.
%

 \subsection{  A priori estimates(I): Time-independent a priori estimates}

In this subsection,
 we will use the convention that $C$ denotes a generic positive constant depending  on the initial data $(\n_0,u_0)$  and some known constants but independent of both $\inf\limits_{x\in\xR}\rho_0(x)$ and $T$, and  use $C(\alpha)$ to emphasize that $C$ depends on $\alpha.$

 We start with the following  energy estimate for the solutions $(\n, u)$.
 \begin{lemma} \la{Rla3.1}
 There is a positive constant $C$ depending  on $\gamma$, $\bar\mu$, $\tilde{\n}$, $\|\n_0-\tilde{\n}\|_{H^1}$,   $\|u_0\|_{H^1}$,  and $ \|\tilde{f}\|_{H^2}$ such that
\be\la{Rgj3.1}
\sup_{0\le t\le T} \xix \left(\n u^2+G(\n)\right)dx+\xiT\xix\mu(\n)u_x^2dxdt\le C,
\ee
where and in what follows $\ti f\triangleq \int_{-\infty}^xf(y)dx$ and $G$ denotes the potential energy density given by
 \be\ba\la{R03.0}
G(\n)&=\n\int_{\tilde{\n}}^\n \frac{P(s)-P(\tilde\n)}{s^2}ds =\frac{1}{\gamma-1}\left[ P(\n)-P(\tilde\n)-\gamma\tilde{\n}^{\gamma-1}(\n-\tilde\n)\right].
\ea\ee
\end{lemma}

\pf  Multiplying $\eqref{R1d}_1$ and $\eqref{R1d}_2$ by  $G'(\n)$ and  $u$ respectively, we obtain after using  integration by parts and  the far-field condition \eqref{R1d2} that
 \bnn\ba\la{R03.1}
& \frac{d}{dt} \xix \left(\frac{1}{2}\n u^2+ G(\n)\right)dx+ \xix\mu(\n)u_x^2dx= \frac{d}{dt}\xix (\n-\tilde \n) \tilde{f}dx,
\ea\enn
which yields
 \be\ba\la{R03.1'}
& \sup_{0\le t\le T}\xix \left(\frac{1}{2}\n u^2+ G(\n) \right)dx+ \xiT\xix\mu(\n)u_x^2dxdt\\
&\le C+ \xix (\n-\tilde \n) \tilde{f} dx- \xix (\n_0-\tilde \n) \tilde{f}dx\\
&\le C+ \xix (\n-\tilde \n) \tilde{f} dx+C\|\n_0-\tilde \n\|_{L^2}^2+C\|\tilde f\|_{L^2}^2\\
&\le C+\int_{R_1} (\n-\tilde \n) \tilde{f} dx+\int_{\xR\setminus R_1} (\n-\tilde \n) \tilde{f} dx, \ea\ee
with
\be\label{gyr1}R_1\triangleq  \{x\in \mathbb{R}~|~\n^{\gamma-1}\le \gamma^2(\tilde \n^{\gamma-1}+1)\}  .\ee

Then, on the one hand, when   $0\le \n\le M$ for some positive constant $M>0$,
there are positive constants $K_1$ and $K_2$ depending only  on $\tilde \n$ and $M$ such that
 \be\ba\la{r03.1}
K_1(\n-\tilde \n)^2 \le G(\n) \le  K_2(\n-\tilde \n)^2,
\ea\ee
which together with Young's inequality implies that
 \be\ba\la{r03.2}
\int_{R_1} |\n-\tilde \n|| \tilde{f} |dx&\le \ve  \|\n-\tilde \n\|_{L^2(R_1)}^2 +C(\ve)\|\tilde f\|_{L^2(R_1)}^2\\&\le C\ve \xix G(\n)dx+C(\ve).
\ea\ee
On the other hand, since $\n>\tilde\n$ in $\xR\setminus R_1,$ we have
 \be\ba\la{r03.3}
\int_{\xR\setminus R_1} |\n-\tilde \n|| \tilde{f}|  dx&\le  2\int_{\xR\setminus R_1} \n |\tilde f| dx\\
&\le \ve \int_{\xR\setminus R_1} \n^{\frac{2\gamma-1}{\gamma}}dx +C(\ve)\int_{\xR\setminus R_1} |\tilde f|^{\frac{2\gamma-1}{\gamma-1}}dx\\
&\le \ve \int_{\xR\setminus R_1} \n^{\frac{2\gamma-1}{\gamma}}dx +C(\ve)\|\tilde f\|_{L^\infty(\xR\setminus R_1)}^{\frac{1}{\gamma-1}} \|\tilde f\|_{L^2(\xR\setminus R_1)}^2\\
&\le \ve \xix G(\n)dx+C(\ve),
\ea\ee
where in the last inequality one has used the following  fact:
 \be\ba\la{r03.5}\notag
\n^{\frac{2\gamma-1}{\gamma}} \le \n\left(\frac{1}{\gamma}\n^{\gamma-1}+\frac{\gamma-1}{\gamma}\right)  \le\frac{1}{\gamma-1}\left(\n^{\gamma}-\gamma\tilde \n^{\gamma-1}\n\right) \le  G(\n),~~~ ~x\in \xR\setminus R_1,
\ea\ee
due to Young's inequality,  \eqref{gyr1}, and \eqref{R03.0}.

Substituting   \eqref{r03.2}  and \eqref{r03.3}  into \eqref{R03.1'}, we obtain      \eqref{Rgj3.1} after  choosing $\ve$  suitably small. This  completes  the proof of Lemma \ref{Rla3.1}.  \hfill $\Box$

Next,
we   derive the key uniform (in time) upper bound of the density,
which is crucial for obtaining both the derivative estimates and the large time behavior of the solutions. This method is motivated by Kazhikhov \cite{ka1}.

 \begin{lemma} \la{Rla3.2}
There is a positive constant $\bar{\n}$ depending  only on $\gamma$, $\bar\mu$, $\tilde{\n}$, $\|\n_0-\tilde{\n}\|_{H^1}$,   $\|u_0\|_{H^1}$,  and $ \|\tilde{f}\|_{H^2}$ such that
  \be\la{Rgj3.2}
0\le \n(x,t)\le \bar{\n},~~~\forall~(x,t)\in  \mathbb{R}\times[0,T].
\ee
\end{lemma}

\pf   For any integer $N,$ let $x\in [N,N+1].$ Integrating \eqref{R1d}$_1$ over $(N,x)$ with respect to  $x $ leads to
\be\la{lj3}\ba  -\mu(\n)u_x+P+\xl(\int_{N}^x\n udy\xr)_t+\n u^2=\xl(-\mu(\n)u_x+P+\n u^2\xr)({N},t)+\int_{N}^x\n fdy,\ea\ee
which in particular implies
\be\la{lj4}\ba &\xl(-\mu(\n)u_x+P+\n u^2\xr)({N},t)\\&= \int_{N}^{{N}+1}\xl(-\mu(\n)u_x+P+\n u^2+\int_{N}^x \n fdy\xr)dx +\xl(\int_{N}^{{N}+1}\int_{N}^x\n udydx\xr)_t\\&\le - \int_{N}^{{N}+1} \mu(\n)u_xdx +C  +\xl(\int_{N}^{{N}+1}\int_{N}^x\n udydx\xr)_t. \ea\ee

due to \eqref{Rgj3.1}  and the following simple fact
\be\la{l9} \ba
  \int_{N}^{{N}+1} \n^\ga dx&\le  C +\int_{N}^{{N}+1}1_{\{\n^{\gamma-1}\ge 2\gamma\tilde\n^{\gamma-1}\}} \n^\ga dy \\
 &\le C +2(\gamma-1)\int_{N}^{{N}+1} G(\n)dy\\
&\le C.
\ea\ee
Denoting
  \bnn
F(\n)\triangleq\int_1^\n \mu(s)s^{-1}ds,
\enn
 one deduces from $\eqref{R1d}_1$ that
  \bnn\la{Rl2}
D_tF(\n)=-\mu(\n)u_x,
\enn
which together with \eqref{lj3} and \eqref{lj4} gives
\be\la{Rl6}\ba D_t(F(\n)+b_1(t))+P&\le  -\int_{N}^{{N}+1}\mu(\n)u_xdx+C \\&\le  \int_{N}^{{N}+1} \mu(\n)dx  \int_{N}^{{N}+1} \mu(\n)u_x^2dx+C\\
 &\le  \sup_{x\in \xR} \frac{\mu(\n)}{\n+1}\int_{N}^{{N}+1} (\n+1)dx \int_{N}^{{N}+1} \mu(\n)u_x^2dx+C\\
&\le  C\sup_{x\in \xR} \frac{\mu(\n)}{\n+1}\int_{-\infty}^\infty \mu(\n)u_x^2dx+C, \ea\ee
where $$b_1(t)\triangleq \int_{N}^x\n udy+\int_{N}^{{N}+1}\int_{N}^x\n udydx$$ satisfies
\be\la{Rl8} \ba
 |b_1(t)| &\le C\int_{N}^{{N}+1}\n| u|dy\\&
 \le C\int_{N}^{{N}+1} \n  dx+C\int_{N}^{{N}+1} \n u^2 dx \le C_1.
\ea\ee due to \eqref{l9} and \eqref{Rgj3.1}.

Since $\mu$ satisfies \eqref{Rn3}, it holds
  \bnn\la{Rn3-1}\ba
\sup_{x\in \xR} \frac{\mu(\n)}{\n+1}&\le  \lambda_0\sup_{x\in \xR}\left(\frac{1}{\n+1}\int_1^\n \mu(s)ds\right)+\lambda_1
\\
 &\le \lambda_0\sup_{x\in \xR}\int_1^{\max\{\n,1\}} \frac{\mu(s)}{s}ds+\lambda_1,
\ea\enn
which together with \eqref{Rl6} leads to
\be\la{Rl4}\ba &D_t(F(\n)+b_1(t))+P \\&\le C_2+C \left(\sup_{x\in \xR}\int_1^{\max\{\n,1\}} \frac{\mu(s)}{s}ds+1\right)\int_{-\infty}^\infty \mu(\n)u_x^2dx,
\ea\ee with   constant $C_2>1$.

Choosing a  constant $\nu\ge C_2^{1/\gamma}$ such that for all $\n\ge \nu$,
\be\la{Rl11} \ba
P(\n)-C_2\ge 0,
\ea\ee and multiplying \eqref{Rl4} by
  \be\la{Rl13} \ba
H\triangleq (F(\n)+ b_1(t)-F(\nu)-C_1)_+,
\ea\ee
we obtain after using \eqref{Rl8}   and \eqref{Rl11} that
  \be\la{Rl14} \ba
D_t H^2\le& CH \left(\sup_{x\in \xR}\int_1^{\max\{\n,1\}} \frac{\mu(s)}{s}ds+1\right)\int_{-\infty}^\infty \mu(\n)u_x^2dx\\
\le &C\sup_{x\in \xR}H^2\xix \mu(\n)u_x^2dx+C \xix \mu(\n)u_x^2dx.
\ea\ee
Integrating \eqref{Rl14} over $(0, t)$ gives that for $x\in [N-1,~N]$,
  \be\la{Rl16} \ba
 H^2 (x ,t)\le \bar C+\bar C\int_0^t \sup_{x\in \xR}H^2 \xix  \mu(\n)u_x^2dxdt,
\ea\ee
where $\bar C$ is a positive constant which is independent of $N$. Since $N$ is arbitrary, the inequality \eqref{Rl16} holds for all $x\in\xR$, that is,
  \bnn
 \sup_{x\in \xR}H^2\le C+C\int_0^t \sup_{x\in \xR}H^2 \xix \mu(\n)u_x^2dxdt,
\enn
which together with Gronwall's inequality and \eqref{Rgj3.1} yields that
  \be\la{Rl18} \ba
 \sup_{(x,t)\in \xR\times(0,\infty)}H^2\le C.
\ea\ee

Consequently, 
  the desired \eqref{Rgj3.2} is a direct consequence of \eqref{Rl13}, \eqref{Rl18}, and  \eqref{Rl8}. This finishes  the proof of Lemma \ref{Rla3.2}. \hfill $\Box$

With the uniform upper bound of the density  \eqref{Rgj3.2} at hand, we have the following time-independent bound  on  $\|\n-\tilde \n\|_{L^2}^2.$
 \begin{corollary}  \la{Rla3.2'} 
 There is a positive constant $C$ depending  on $\gamma$, $\bar\mu$, $\tilde{\n}$, $\|\n_0-\tilde{\n}\|_{H^1}$,   $\|u_0\|_{H^1}$,  and $ \|\tilde{f}\|_{H^2}$ such that
  \be\la{Rgj3.2'}
\sup_{0\le t\le T}\left( \|\n-\tilde \n\|_{L^2}^2+\|P(\n)-P(\tilde \n)\|_{L^2}^2 \right)\le C.
\ee
\end{corollary}

\pf  It follows from \eqref{Rgj3.1} and  \eqref{r03.1}  that
\bnn
\ba\la{R3.21}
   &\|\n-\tilde \n\|_{L^2}^2+\|P(\n)-P(\tilde \n)\|_{L^2}^2\\
   &\le C  \xix G(\n)dx+\|(\gamma-1)G(\n) +\gamma\tilde\n^{\gamma-1}\n-\gamma\tilde\n^{\gamma}\|_{L^2}^2\\
   &\le C \xix G(\n)dx \le C,
\ea
\enn
which completes the proof of Corollary \ref{Rla3.2'}.  \hfill $\Box$

 \subsection{\la{Rse3} A priori estimates(II): Time-dependent a priori estimates}
 In this subsection, we now proceed to derive  the  derivative estimates of
 the strong  $(\n,u)$ to the Cauchy problem \eqref{R1d}--\eqref{R1d2}. In what follows,  the constant $C$   depends on $T$ but is still independent of $\inf\limits_{x\in \xR}\rho_0(x).$


\begin{lemma} \la{Rla3.5} 
 There is a positive constant $C$ depending  on $T$, $\gamma$, $\bar\mu$, $\tilde{\n}$, $\|\n_0-\tilde{\n}\|_{H^1}$,   $\|u_0\|_{H^1}$,  and $ \|\tilde{f}\|_{H^2}$ such that
 \be
\ba\la{Rgj3.5}
   \sup_{0\le t\le T}   \| u \|_{H^1}^2 +\xiT \|\n^{1/2}\dot u\|_{L^2}^2 dt\le C(T).
\ea
\ee
\end{lemma}

\pf
First, multiplying $\eqref{R1d}_{2}$  by  $ \dot u$ and  integrating the resulting equation by parts yield
\be
\ba\la{R3.2}
  &\frac{1}{2}\frac{d}{dt}\xix \mu(\n)u_x^2dx+\xix\n\dot u^2dx\\
  &=\frac{d}{dt}\left(\xix  [P(\n)-P(\tilde\n)] u_xdx+\xix  \n fudx\right)\\
  &\quad-\frac{1}{2}\xix [\mu(\n)+\mu'(\n)\n]u_x^3dx+\ga\xix P(\n)u_x^2dx-\xix \n u^2 f_xdx\\
  &\le\frac{d}{dt}\left(\xix [P(\n)-P(\tilde\n)] u_xdx+\xix \n  fudx\right)\\
  &\quad+C \|u_x\|_{L^3}^3+C  \|u_x\|_{L^2}^2+ C\|u\|_{L^\infty}\|\n^{1/2}u\|_{L^2}\|f_x\|_{L^2} \\
    &\le\frac{d}{dt}\left(\xix[P(\n)-P(\tilde\n)]u_xdx+\xix \n fudx\right) \\
    &\quad+C \|u_x\|_{L^\infty}\|u_x\|_{L^2}^2+C\|u_x\|_{L^2}^2+C,
\ea
\ee
where in the last inequality one has used the following fact:
\be
\ba\la{R3.5}
   \|u\|_{L^2}^2+ \|u\|_{L^\infty}
   =&\tilde\n^{-1}\left(\xix \n u^2dx+\xix (\tilde\n-\n)u^2dx\right)+ \|u\|_{L^\infty}\\
   \le &C\|\n^{1/2}u\|_{L^2}^2+C\|\n-\tilde\n\|_{L^2}\|u\|_{L^2}\|u\|_{L^\infty} + \|u\|_{L^\infty}\\
   \le &C\|\n^{1/2}u\|_{L^2}^2+C\|u\|_{L^2}^{3/2} \|u_x\|_{L^2}^{1/2}+C\|u\|_{L^2}^{1/2} \|u_x\|_{L^2}^{1/2}\\
   \le &C+C\|u_x\|_{L^2}^2+\frac{1}{2}\|u\|_{L^2}^2.
\ea
\ee
due to \eqref{Rgj3.2'},  \eqref{Rgj3.1},
and the Sobolev inequality.
Furthermore, using \eqref{Rn3}, \eqref{Rgj3.1}, \eqref{R1d}$_2$,  and \eqref{Rgj3.2'}, we  get
 \be
\ba\la{R3.7}
  \|u_x\|_{L^\infty}
  \le & C\|\mu(\n)u_x-P(\n)+P(\tilde\n)\|_{L^\infty}+ C\|P(\n)-P(\tilde\n)\|_{L^\infty}\\
    \le &C  \|\mu(\n)u_x-[P(\n)-P(\tilde\n)]\|_{L^2}+C  \|(\mu(\n)u_x- P(\n))_x  \|_{L^2}+C \\
    \le & C\|\mu(\n)u_x\|_{L^2}+C\|P(\n)-P(\tilde\n)\|_{L^2}+C\|\n\dot u\|_{L^2}+C\|\n f\|_{L^2}+C \\
    \le & C \|\sqrt{\mu(\n)}u_x\|_{L^2}+C \|\n^{1/2}\dot u\|_{L^2}+C .
\ea
\ee

Then, putting  \eqref{R3.7} into \eqref{R3.2}, we obtain after   using Young's inequality that
 \be
\ba\la{R3.8}
  & \frac{d}{dt} B(t)+\frac{1}{2}\xix\n\dot u^2dx\le C+C\|\sqrt{\mu(\n)}u_x\|_{L^2}^2+C\|\sqrt{\mu(\n)}u_x\|_{L^2}^4,
\ea
\ee
where
\be \label{rjia1}
 B(t)\triangleq \frac{1}{2}\xix \mu(\n)u_x^2dx-\xix[P(\n)-P(\tilde\n)]u_xdx-\xix \n fudx
\ee
satisfies
 \be
\ba\la{R3.6}
\frac{1}{4}\|\sqrt{\mu(\n)}u_x\|_{L^2}^2-C\le  B(t)\le  C\|\sqrt{\mu(\n)}u_x\|_{L^2}^2+C.
\ea
\ee
Indeed, it follows from  \eqref{Rn3}, \eqref{Rgj3.1}, and \eqref{Rgj3.2'} that \bnn
\ba\la{R3.10}
  & \xix [P(\n)-P(\tilde\n)]  u_xdx+ \xix \n fudx\\
     &\le \frac{1}{4}\xix \mu(\n)u_x^2dx+C\|P(\n)-P(\tilde\n)\|_{L^2}^2+C \|\n^{1/2}u\|_{L^2}\|f\|_{L^2} \\
     &\le \frac{1}{4}\xix \mu(\n)u_x^2dx+C,
\ea
\enn which along with \eqref{rjia1} yields \eqref{R3.6}.

Finally,  Gronwall's inequality combined with \eqref{R3.8}, \eqref{R3.6}, \eqref{Rn3},  \eqref{Rgj3.1}, and \eqref{R3.5} yields \eqref{Rgj3.5} and completes the proof of Lemma \ref{Rla3.5}.  \hfill $\Box$
%

\begin{lemma} \la{Rla4.1}
 There is a positive constant $C$ depending  on $T$, $\gamma$, $\bar\mu$, $\tilde{\n}$, $\|\n_0-\tilde{\n}\|_{H^1}$,   $\|u_0\|_{H^1}$,  and $ \|\tilde{f}\|_{H^2}$ such that
 \be
\ba\la{Rgj4.1}
   \sup_{0\le t\le T} \left(\|\n_x\|_{L^2}^2+\|\n_t\|_{L^2}^2\right) \le C ,
\ea
\ee
\end{lemma}

\pf Differentiating $\eqref{R1d}_1$ with respect to $x$ gives
\be
\ba\la{Ra1}
   \n_{xt}+\n_{xx}u+2\n_xu_x+\n u_{xx}=0.
\ea
\ee
Multiplying \eqref{Ra1} by $\n_x$  and integrating the resultant equation by parts  yield  that
\be\ba\la{Ra2}
 \frac{d}{dt}\|\n_x\|_{L^2} 
   &\le C\|u_x\|_{L^\infty}\|\n_x\|_{L^2} +C \|u_{xx}\|_{L^2},
\ea
\ee
due to \eqref{Rgj3.2}.

Then, it follows from $\eqref{R1d}_2$ that
\bnn
\ba\la{Ra7}
\mu(\n)u_{xx}= \n \dot u+P_x-\n f-\mu'(\n)\n_xu_x,
\ea
\enn
which together with \eqref{Rn3},  and  \eqref{Rgj3.2} yields that
\be
\ba\la{Ra9}
\|u_{xx}\|_{L^2}&\le  C \| \n^{1/2} \dot u\|_{L^2}+ C \| \n_x\|_{L^2}+ C \|f\|_{L^2}+C \| \n_x\|_{L^2}\| u_x\|_{L^\infty} .
\ea
\ee
Submitting \eqref{Ra9} into \eqref{Ra2}, one obtains after using  \eqref{R3.7},  \eqref{Rgj3.5}, and Gronwall's inequality
\be
\ba\la{Ra5}
   \sup_{0\le t\le T}  \|\n_x\|_{L^2}  \le C ,
\ea
\ee
which along with $\eqref{R1d}_1$, \eqref{Rgj3.2}, \eqref{R3.5},    and \eqref{Rgj3.5}  leads to
\bnn
\ba\la{Ra6}
 \|\n_t\|_{L^2}  & \le C\|u\|_{L^\infty} \|\n_x\|_{L^2} +C \| u_x\|_{L^2} \le C .
\ea
\enn
Combining this  with \eqref{Ra5} shows  \eqref{Rgj4.1} and
finishes the proof of Lemma \ref{Rla4.1}.  \hfill $\Box$

\begin{lemma} \la{Rla4.1'}There is a positive constant $C$ depending  on $T$, $\gamma$, $\bar\mu$, $\tilde{\n}$, $\|\n_0-\tilde{\n}\|_{H^1}$,   $\|u_0\|_{H^1}$,  and $ \|\tilde{f}\|_{H^2}$ such that
\be
\ba\la{nRgj4.1}
   \sup_{0\le t\le T}  \sigma\| u_{xx}\|_{L^2}^2 +\int_0^T\left(\| u_{xx}\|_{L^2}^2+\si\|u_{t}\|_{H^1}^2  \right)dt\le C,
\ea
\ee with  $\si(t)\triangleq\min \{1,t\}.$
\end{lemma}

\pf  First, operating $\pa/\pa_t+(u\cdot~)_x$ to $ (\ref{R1d})_2 $ yields
 \be\ba\la{R3.12}
   \n \dot u_t+\n u\dot u_x-[\mu(\n) \dot u_x]_x=-\gamma [P(\n)u_x]_x-[(\mu(\n)+\mu'(\n)\n)u_x^2]_x+\n uf_x.
 \ea\ee
 Multiplying \eqref{R3.12} by $\dot{u} $, one gets after using \eqref{R3.5}, \eqref{R3.7}, and \eqref{Rgj3.5}   that
 \be
\ba\la{R3.15}
  &\frac{1}{2}\frac{d}{dt}\xix \n|\dot u|^2dx+\xix \mu(\n) |\dot u_x|^2dx\\
  &=\gamma \xix P(\n)u_x\dot u_xdx+\xix (\mu(\n)+\mu'(\n)\n)u_x^2 \dot u_xdx +\xix \n uf_x \dot udx\\&\le C\|\dot u_x\|_{L^2}+ C\|\dot u_x\|_{L^2}\|u_x\|_{L^\infty}+C \|u\|_{L^\infty}\|f_x\|_{L^2}\|\n^{1/2}\dot u\|_{L^2}\\&\le \frac12\|\sqrt{\mu(\n)}\dot u_x\|^2_{L^2}+ C \|\n^{1/2}\dot u\|_{L^2}^2+C .
\ea
\ee

Then, multiplying \eqref{R3.15} by $\sigma$ and integrating the resulting inequality  over $(0,T)$,   we obtain after using \eqref{Rgj3.5}  that
 \be
\ba\la{R3.16}
   &\sup_{0\le t\le T} \sigma\xix \n|\dot u|^2dx+\xiT\sigma\xix \mu(\n) |\dot u_x|^2dxdt \le C(T),
\ea
\ee
which along with  \eqref{R3.5} and  \eqref{R3.7} implies that
 \be
\ba\la{Rgj3}
   &\sup_{0\le t\le T} \sigma\|u_x\|_{L^\infty}^2 \le C(T).
\ea
\ee

 Finally, it follows from direct calculations, \eqref{Rgj3.5}, and \eqref{R3.5}   that
\bnn
\ba\la{Rb5}
 \|\n^{1/2}u_t\|_{L^2}^2
 &\le C\|\n^{1/2}\dot u\|_{L^2}^2+  C\|\n\|_{L^\infty} \|  u\|_{L^\infty}^2 \|u_x \|_{L^2}^2\\&\le C \|\n^{1/2}\dot u\|_{L^2}^2+ C,\ea\enn which gives \bnn\ba
 \|u_{ t}\|_{H^1}^2&\le   C\int |\ti\n-\n|u_t^2dx+C \int \n u_t^2dx+C\|u_{xt}\|_{L^2}^2 \\&\le C\|\n-\ti\n\|_{L^2}\|u_{ t}\|_{L^2}^{3/2}\|u_{ xt}\|_{L^2}^{1/2}+C\|\n^{1/2}\dot u\|_{L^2}^2+ C+C\|u_{xt}\|_{L^2}^2\\&\le \frac12 \|u_{ t}\|_{L^2}^2+C\|\n^{1/2}\dot u\|_{L^2}^2+ C+C\|u_{xt}\|_{L^2}^2 \\ &=   \frac12 \|u_{ t}\|_{L^2}^2+C\|\n^{1/2}\dot u\|_{L^2}^2+ C+\|(\dot u-uu_x)_x\|_{L^2}^2 \\
 &\le \frac12 \|u_{ t}\|_{L^2}^2+C\|\n^{1/2}\dot u\|_{L^2}^2+ C+C\|\dot u_x\|_{L^2}^2\\
 &\quad + \|u_x\|_{L^\infty}^2\|u_x\|_{L^2}^2+ \|u\|_{L^\infty}^2\|u_{xx}\|_{L^2}^2\\
 &\le \frac12 \|u_{ t}\|_{L^2}^2+C\|\dot u_x\|_{L^2}^2+C\|\n^{1/2}\dot u\|_{L^2}^2+C,
\ea
\enn where in the last inequality we have used \eqref{R3.5}, \eqref{Rgj3.5}, \eqref{Ra9}, \eqref{Rgj4.1}, and \eqref{R3.7}. Combining this, \eqref{Ra9}, \eqref{R3.5}, \eqref{R3.16},  and \eqref{Rgj3} gives   \eqref{nRgj4.1}.
The proof of Lemma \ref{Rla4.1} is completed. \hfill $\Box$

\section{Proof of Theorem  \ref{Rth1} }

 First, for $(\n_0,u_0,\tilde{f})$   satisfying the conditions as  in Theorem \ref{Rth1},
we construct  the smooth approximate data as follows:
\be\la{8.1}
\n_0^{\delta,\eta}=\frac{\n_0\ast j_\delta+\eta\tilde \n}{1+\eta},~~~ u_0^{\delta,\eta}=u_0\ast j_\delta,~~~ \tilde{f}^{\delta,\eta}=\tilde{f}\ast j_\delta,\ee
where $ \delta\in (0,1),\eta\in (0,1),$ and $j_\delta(x)$ is the standard mollifier with width $\delta$. It is easy to check that
\be\la{8.2}\begin{cases}
0<\frac{\eta\tilde \n}{1+\eta}\le\n_0^{\delta,\eta}\le \frac{\sup\limits_{x\in\xR}\n_0 +\eta\tilde \n}{1+\eta}<\infty,\\
\lim_{\delta,\eta\rightarrow 0}\left(\|\n_0^{\delta,\eta}-\n_0\|_{H^1}+\|u_0^{\delta,\eta}-u_0\|_{L^2}+\|\tilde{f}^{\delta,\eta}-\tilde{f}\|_{H^2}\right)=0.
\end{cases}\ee

Choosing  $   \mu_{\eta}\in C^2[0,\infty) $ satisfying $\lim_{\eta\rightarrow 0}\|\mu_{\eta}-\mu\|_{C^1[0,M]}$ for any $M>0,$  we consider  the Cauchy problem \eqref{R1d}-\eqref{R1d2}  with $\mu$ replaced by $\mu_\eta$ and the data $(\n_0^{\delta,\eta},~u_0^{\delta,\eta},~\tilde{f}^{\delta,\eta})$ satisfying \eqref{8.1}--\eqref{8.2}.  It follows  from Lemma \ref{Rpro1} that there exists a unique strong solution $(\n , u )$  of problem \eqref{R1d}--\eqref{R1d2} on $\xR\times[0,T_{\de,\eta}]$.  Moreover,  the estimates obtained in Lemmas \ref{Rla3.1}--\ref{Rla4.1}  show that the solution
$(\n , u )$ satisfies for any $0<T\le T_{\de,\eta}$,
\be\ba
   \la{R8.6}
&  \sup_{0\le t\le T} \xl(\|(\n -\ti\n,\mu(\n )-\mu(\ti\n),P(\n )-P(\ti\n))\|_{H^1}
+\| \n _t\|_{L^2}\right.\\
&\left.+\| \n  u  \|_{L^2}+\|u \|_{H^1}+\sqrt{t}\|\sqrt{\n }u _t\|_{L^2}
+\sqrt{t}\|u _{xx}\|_{L^2}\xr)\\& +\xiT\xl( \|u \|_{H^2}^2+ t\|u _{xt}\|_{L^2}^2\xr)dt\le \bar C, \ea\ee
where $\bar C$ is independent of $\delta$ and $\eta$. Moreover, similar to   \cite{lww}, we can prove that there exists some $\tilde C$ depending on $\de$ and $\eta$ such that  \bnn \sup_{0\le t\le T} \left( \|\n_{xx}\|_{L^2}^2 + \|P_{xx}\|_{L^2}^2 +\|\n_{xt}\|_{L^2}^2 + \|P_{xt}\|_{L^2}^2 +t\|u_{xt}\|_{L^2}^2 + t\|u_{xxx}\|_{L^2}^2 \right)\le \tilde C,\enn which in particular implies that  \be
   \la{1d5''}\begin{cases}
   \rho-\ti\n, ~P(\n )-P(\ti\n)  \in C([0,T]; H^2), \quad
    u\in C([0,T]; H^2) \cap L^2(0,T;H^3), \\
 u_t\in L^2(0,T; H^1),~~t^{1/2}u \in L^\infty(0,T; H^3),\\
t^{1/2}u_t \in L^\infty(0,T; H^1)\cap L^2(0,T; H^2),~~  t^{1/2}\n^{1/2}u_{tt}\in L^2(0,T; L^2).\end{cases} \ee

Then, we will extend the   existence time $T_{\de,\eta}$ to be infinity. Indeed, let $T^*$ be the  maximal time of existence for the strong solution. Then, $T^*\ge T_{\de,\eta}$.    If $T^*<\infty,$
defining
  \bnn \la{nx3} (\n^*, u^*)\triangleq(\n, u)(x,T^*)=\lim_{t\rightarrow T^*}(\n, u)(x,t),\enn
we can derive from \eqref{1d5''}   that $ (\n^*,P^*,u^*)$ satisfies the initial condition \eqref{R1d3} at $t=T^*$.   Therefore, one can take $ (\n^*, u^*)$ as the initial data at $t=T^*$ and then use the local existence theory (Lemma \ref{Rpro1}) to extend  the strong solution beyond the maximum existence time $T^*$. This contradicts the assumption on $T^*$. Thus, $T^*=\infty.$

Finally, we rewrite the global strong solutions  on $\xR\times [0,\infty)$  obtained above as $(\n^{\delta,\eta},u^{\delta,\eta}) .$ With the estimates \eqref{R8.6} at hand,  letting first $\delta\rightarrow 0$ and then $\eta\rightarrow 0$, we  find  that the sequence
$(\n^{\delta,\eta},u^{\delta,\eta})$ converges, up to the extraction of subsequences, to some limit $(\n,u)$   in the
obvious weak sense.  Then we deduce from \eqref{R8.6} that  $(\n,u)$ is a strong  solution of \eqref{R1d}-\eqref{R1d2}   on $\mathbb{R}\times (0,T]$ (for any $0<T<\infty$) satisfying \eqref{R1d5}.
   Moreover, the uniqueness of the strong solution $(\n,u)$ is guaranteed by the regularities \eqref{R1d5}. For the detailed proof,  please see \cite{liliang}.
The proof of Theorem  \ref{Rth1} is completed. \hfill $\Box$

\section{Proof of Theorem \ref{Rth2}}

  We will divide the proof into three steps.

\textbf{Step 1.} For the case of  $f\equiv 0$ in $\eqref{R1d}$ , we will derive some time-independent lower order estimates of the solutions $(\n,u)$ obtained in Lemmas \ref{Rla3.5}--\ref{Rla4.1'}, that is
 \be
\ba\la{xRgj3.5}
   \sup_{0\le t\le T} \left(\|u\|_{L^2\cap L^\infty}+ \|\sqrt{\mu(\n)} u_x\|_{L^2}^2\xr)+\xiT \|\n^{1/2}\dot u\|_{L^2}^2 dt\le C,
\ea
\ee \be
\ba\la{xRgj3.6}
   \sup_{0\le t\le T}   \sigma\left( \|u_x\|_{L^\infty}^2+\|\n^{1/2}\dot u\|_{L^2}^2\right)+\xiT\sigma\| \sqrt{\mu }\dot u_x\|_{L^2}^2 dt\le C.
\ea
\ee

First,  multiplying $\eqref{R1d}_{2}$ with $f\equiv 0$  by  $ \dot u$ and  integrating   by parts, it follows from the same arguments as the proof of \eqref{R3.2} that
\be
\ba\la{xR3.2}
  &\frac{1}{2}\frac{d}{dt}\xix \mu(\n)u_x^2dx+\xix\n\dot u^2dx\\
  &\le\frac{d}{dt}\xix [P(\n)-P(\tilde\n)] u_xdx+C \|u_x\|_{L^\infty}\|u_x\|_{L^2}^2+C\|u_x\|_{L^2}^2\\
  &\le  \frac{d}{dt}\xix [P(\n)-P(\tilde\n)] u_xdx+\frac{1}{2}\|\n^{1/2}\dot u\|_{L^2}^2+C \|u_x\|_{L^2}^4+C\|u_x\|_{L^2}^2,
\ea
\ee
where in the last inequality one has used \eqref{R3.7} and Young's inequality. Since
\be\la{rgy9}\xix [P(\n)-P(\tilde\n)] u_xdx\le \frac{\bar \mu}{4}\|u_x\|_{L^2}^2+\|P(\n)-P(\tilde\n)\|_{L^2}^2\le \frac{\bar \mu}{4}\|u_x\|_{L^2}^2+C\ee
due to \eqref{Rgj3.2'}, the Gronwall's inequality together with \eqref{rgy9},  \eqref{xR3.2}, and \eqref{Rgj3.1} yields \eqref{xRgj3.5}.

Next, operating $(\pa/\pa_t+(u\cdot~)_x)\dot u$ to \eqref{R1d} with $ f\equiv 0 $  and integrating the resulting equation by parts, we deduce from the same calculations as \eqref{R3.12}--\eqref{R3.15} that
 \be
\ba\la{xR3.13}
  &\frac{1}{2}\frac{d}{dt}\xix \n|\dot u|^2dx+\xix \mu(\n) |\dot u_x|^2dx\\
  &=\gamma \xix P(\n)u_x\dot u_xdx+\xix (\mu(\n)+\mu'(\n)\n)u_x^2 \dot u_xdx  \\
  &\le   \frac{\bar \mu}{2}\|\dot u_x\|_{L^2}^2 +C\|u_x\|_{L^2}^2\|\n^{1/2}\dot u\|_{L^2}^2+C\|u_x\|_{L^2}^4+C\|u_x\|_{L^2}^2.
\ea
\ee
Multiplying \eqref{xR3.13}  by  $\si$, one gets after using  Gronwall's inequality, \eqref{R1d4}, \eqref{Rgj3.1}, and \eqref{xRgj3.5} that
 \be
\ba\la{xR3.16}
  \sup_{0\le t\le T} \si\|\n^{1/2}\dot u\|_{L^2}^2+\xiT \si\|\dot u_x\|_{L^2}^2dt
  &\le C,
\ea
\ee
which together with  \eqref{R3.5},  \eqref{R3.7},  and \eqref{xRgj3.5} implies  that
 \be
\ba\la{xRgj3}
   &\sup_{0\le t\le T} \left(\|u\|_{L^2\cap L^\infty}+\si\|u_x\|_{L^\infty}\right)\le C.
\ea
\ee
The combination of  \eqref{xRgj3} with  \eqref{xR3.16} gives  \eqref{xRgj3.6}.

\textbf{Step 2.}   We will show that for any $p\ge 2$,
\be\la{Rgj3.4'}
\lim_{t\rightarrow0} \|u_x\|_{L^p} =0.
\ee

First, we claim that
\be\label{Rts1}
\lim_{t\rightarrow\infty}\int_{t-1}^t\left(\|u_x\|_{L^2}^2+\left|\frac{d}{d\tau}\|u_x\|_{L^2}^2\right|\right)d\tau=0.
\ee

For any $t>1$ and $s\in(t-1,t)$, it holds
\bnn
\|u_x\|_{L^2}^2(t)-\|u_x\|_{L^2}^2(s)=\int_t^s \frac{d}{d\tau}\|u_x\|_{L^2}^2d\tau\le \int_{t-1}^t  \left|\frac{d}{d\tau}\|u_x\|_{L^2}^2\right| d\tau,
\enn
that is
\bnn\label{Rts2}
\|u_x\|_{L^2}^2(t)\le \int_{t-1}^t  \left(\|u_x\|_{L^2}^2 +\left|\frac{d}{d\tau}\|u_x\|_{L^2}^2\right|\right) d\tau,
\enn
which together with \eqref{Rts1} yields
\be\label{Rts3}
\lim_{t\rightarrow\infty}\|u_x\|_{L^2}^2(t)=0.
\ee
It is easy to deduce from \eqref{xRgj3}  that  for any  $p>2$,
\be\ba\label{Rts6}
 \|u_x\|_{L^p}^p\le  \|u_x\|_{L^2}^2\|u_x\|_{L^\infty}^{p-2} \le  C\|u_x\|_{L^2}^2,~~~~\forall~~t>1.
\ea\ee
 The combination of \eqref{Rts3} with \eqref{Rts6} yields \eqref{Rgj3.4'}.

Now, it remains  to prove \eqref{Rts1}. Indeed, it follows from integration by parts and \eqref{xRgj3} that
\bnn\ba\label{Rts7}
\left| \frac{d}{dt}\|u_x\|_{L^2}^2\right|& = \left|2\xix u_xu_{xt}dx\right|=\left|2\xix u_x\left(\dot u_x-(uu_x)_x\right)dx \right|\\
 &= \left|2\xix u_x \dot u_xdx- \xix u_x^3 dx\right|\\
 &\le  C\|\dot u_x\|_{L^2}^2+C(1+\|u_x\|_{L^\infty})\|u_x\|_{L^2}^2\\
  &\le  C\|\dot u_x\|_{L^2}^2+C\|u_x\|_{L^2}^2.
\ea\enn
   Combining this, \eqref{Rgj3.1}, and \eqref{xRgj3.6} yields
\bnn\ba\label{Rts8}
\int_1^\infty\left(\|u_x\|_{L^2}^2+\left|\frac{d}{dt}\|u_x\|_{L^2}^2\right|\right)dt\le C,
\ea\enn
which yields  \eqref{Rts1}. The proof of step 2 is completed.

\textbf{Step 3.} Using the methods due to \cite{hlx1,lzz}, we are now in a position to prove that for any $p>2$,
\be\ba\la{Rgj3.4}   \lim_{t\rightarrow\infty}\|\n-\tilde \n\|_{L^p}=0.\ea\ee

First, we claim that
\be\ba\label{Rst21}
\int_0^\infty  \|\n-\tilde \n\|_{L^6}^6 dt\le C.\ea\ee
It follows from integration by parts, \eqref{R1d}$_1,$ and  \eqref{Rgj3.2}  that
\bnn\ba\label{Rst22}
 \frac{d}{dt}\|\n-\tilde \n\|_{L^6}^6  &=6\xix(\n-\tilde \n)^5(\n-\tilde \n)_tdx\\
 &=-6\xix(\n-\tilde \n)^5(\n-\tilde \n)_xudx-6\xix(\n-\tilde \n)^5\n u_xdx\\
 &\le C \|\n-\tilde \n\|_{L^6}^3\|u_x\|_{L^2}+C\|\n-\tilde \n\|_{L^6}^3\|u_x\|_{L^2}\\
 &\le C\|\n-\tilde \n\|_{L^6}^6+C\|u_x\|_{L^2}^2,\ea\enn
which together with \eqref{Rst21} and \eqref{Rgj3.1} gives
\bnn\ba\label{Rst23}
\int_0^\infty  \left|\frac{d}{dt}\|\n-\tilde \n\|_{L^6}^6 \right| dt\le C.\ea\enn
Combining this   with \eqref{Rst21} implies that
\bnn\ba\label{Rst24}   \lim_{t\rightarrow\infty}\|\n-\tilde \n\|_{L^6}=0, \ea\enn
which along with \eqref{Rgj3.2'}  and  \eqref{Rgj3.2} leads to the desired \eqref{Rgj3.4}.

Next, we need   to prove \eqref{Rst21}.
For
\be\ba\label{Rst3}
Q\triangleq \mu(\n)u_x-(P(\n)-P(\tilde \n)),\ea\ee
it   follows from \eqref{R1d}$_2$ with $f\equiv 0$ that
\bnn
Q_{x }= \n\dot u ,\enn  which gives
\be\ba\label{Rst5}
\|Q_x\|_{L^2}=\|\n\dot u\|_{L^2}\le C\|\n^{1/2}\dot u\|_{L^2}.\ea\ee
It follows from \eqref{R1d}$_1$ that
\be\ba\label{Rst1}
(\n-\tilde \n)_t+(\n-\tilde \n)_xu+(\n-\tilde \n)u_x+\tilde \n u_x=0.\ea\ee
Multiplying \eqref{Rst1} by $6(\n-\tilde \n)^5$ and integrating by parts implies
\be\ba\label{Rst2}
&\frac{d}{dt}\|\n-\tilde \n\|_{L^6}^6 +
6\xix \tilde \n(\mu(\n))^{-1}(\n-\tilde \n)^5 (P(\n)-P(\tilde \n))dx\\
&=-
6\xix \tilde \n(\mu(\n))^{-1}(\n-\tilde \n)^5 Q dx-5\xix (\n-\tilde \n)^6u_xdx\\
&\le C \|\n-\tilde \n\|_{L^6}^{5}\|Q \|_{L^6}+C\|\n-\tilde \n\|_{L^\infty}^3\|\n-\tilde \n\|_{L^6}^3\|u_x\|_{L^2}\\
&\le \ve \|\n-\tilde \n\|_{L^6}^6+C( \ve)\|Q\|_{L^2}^4\|Q_x\|_{L^2}^2+C( \ve)\|u_x\|_{L^2}^2
\\
&\le \ve \|\n-\tilde \n\|_{L^6}^6+C( \ve)\|\n^{1/2}\dot u\|_{L^2}^2+C( \ve)\|u_x\|_{L^2}^2,\ea\ee
where one has used \eqref{Rst3}, \eqref{Rn3}, \eqref{Rgj3.2}, \eqref{Rst5}, \eqref{xRgj3.6},  and \eqref{Rgj3.2'}. Furthermore, the direct  calculations combined with   \eqref{Rn3} and \eqref{Rgj3.2} show that for some $0<\alpha<1$,
\be\ba\label{Rst6}
6 \tilde \n(\mu(\n))^{-1}(\n-\tilde \n)^5 (P(\n)-P(\tilde \n))&=6\tilde \n(\mu(\n))^{-1}(\n-\tilde \n)^6 P'(\alpha\n+(1-\alpha)\tilde\n) \\
 &\ge   C_0(\n-\tilde \n)^6,\ea\ee
where the positive constant $ C_0$ depending only on $\gamma, \bar\mu, \bar\n$ and $\tilde\n$.
Substituting \eqref{Rst6} into \eqref{Rst2} and choosing $\ve$ suitably small yield
\be\ba\label{Rst7}
 \frac{d}{dt}\|\n-\tilde \n\|_{L^6}^6 + C_0\|\n-\tilde \n\|_{L^6}^6\le C\|\n^{1/2}\dot u\|_{L^2}^2+C\|u_x\|_{L^2}^2.\ea\ee
 Thus, one can derive the desired \eqref{Rst21} from \eqref{Rst7}, \eqref{Rgj3.1}, \eqref{Rgj3.2'},  and \eqref{xRgj3.5}.

Finally, noticing that \bnn \|\n-\ti\n\|_{C(\overline{\xR})}\le C\|\n-\ti\n\|_{L^6}^{3/4}\|\n_x\|_{L^2}^{1/4},\enn the proof of \eqref{bp} is similar as that of \cite[Theorem 1.2]{lijde}(see also \cite{hlx1}).
 The proof of Theorem \ref{Rth2} is finished. \hfill$\Box$

\end{document}